\documentclass[10pt]{amsart}
\usepackage{amssymb}
\usepackage[OT2,OT1]{fontenc}

\let\newpf\proof \let\proof\relax 
\newenvironment{pf}{\newpf[\proofname]}{\qed\endtrivlist}

\def\DC{{\mathrm{DC}}}

\def\bm{\begin{matrix}}
\def\em{\end{matrix}}

\newcommand{\bt}{\begin{thm}}
\newcommand{\et}{\end{thm}}

\newcommand{\bl}{\begin{lemma}}
\newcommand{\el}{\end{lemma}}

\def\Diff{{\mathrm{Diff^\infty(\R/\Z)}}}

\newcommand{\beq}{\begin{eqnarray}}
\newcommand{\eeq}{\end{eqnarray}}

\def\be{\begin{equation}}
\def\ee{\end{equation}}

\def\ba{{\begin{align}}}
\def\ea{{\end{align}}}

\def\0{{\mathbf 0}}

\newtheorem{thm}{Theorem}[section]

\newtheorem{lemma}[thm]{Lemma}

\theoremstyle{remark}

\numberwithin{equation}{section}

\def \bn {\hfill \\ \smallskip\noindent}

\theoremstyle{definition}

\newtheorem{definition}{Definition}[section]
\def\proof{\bn {\bf Proof.} }

\def\note#1
{\marginpar
{\tiny $\leftarrow$
\par
\hfuzz=20pt \hbadness=9000 \hyphenpenalty=-100 \exhyphenpenalty=-100
\pretolerance=-1 \tolerance=9999 \doublehyphendemerits=-100000
\finalhyphendemerits=-100000 \baselineskip=6pt
#1}\hfuzz=1pt}

\newcommand{\dist}{\operatorname{dist}}

\newcommand{\inter}{\operatorname{int}}

\newcommand{\id}{\operatorname{id}}

\newcommand{\N}{{\mathbb N}}

\newcommand{\R}{{\mathbb R}}

\newcommand{\Z}{{\mathbb Z}}

\def\B0{{\bold{0}}}

\catcode`\@=12

\def\Empty{}
\newcommand\oplabel[1]{
  \def\OpArg{#1} \ifx \OpArg\Empty {} \else
  	\label{#1}
  \fi}
		
\newcommand{\comm}[1]{}
\newcommand{\comment}[1]{}

\begin{document}

\bigskip\bigskip

\title{Distortion elements in $\Diff$}

\author {Artur Avila}

\address{
CNRS UMR 7599, Laboratoire de Probabilit\'es et Mod\`eles al\'eatoires\\
Universit\'e Pierre et Marie Curie--Boite courrier 188\\
75252--Paris Cedex 05, France
}
\curraddr{IMPA, Estrada Dona Castorina 110, Rio de Janeiro, 22460-320,
Brazil}
\email{artur@math.sunysb.edu}

\date{\today}

\begin{abstract}
We consider the group of smooth diffeomorphisms of the circle.
We show that
any recurrent $f$ (in the sense that $\{f^n\}_{n \in \Z}$ is not discrete)
is in fact a distortion element (in the sense that its
iterates can be written as short compositions involving finitely many
smooth diffeomorphisms).  Thus rotations are distortion elements.
\end{abstract}

\setcounter{tocdepth}{1}

\maketitle

\section{Introduction}

Given a finitely generated group $G$ and a choice of a generating
set $S \subset G$, we define a metric $d_S$ on $G$ by taking
$d_S(g_1,g_2)$ as the infimum over all $k \geq 0$ such that there exist
$s_1,...,s_k \in S$ and $\varepsilon_1,...,\varepsilon_k \in \{1,-1\}$ such
that $g_2=s_1^{\varepsilon_1} \cdots s_k^{\varepsilon_k} \cdot g_1$.
Different choices of a generating set give
rise to quasi-isometric metrics.

\begin{definition}

Let $G$ be a group.  We say that $f \in G$ is a
{\it distortion element} if there exists a finitely generated subgroup
$G' \subset G$ containing $f$ such that $d_S(f^n,\id)=o(n)$ for some
(and automatically for every) generating set $S \subset G'$.

\end{definition}

While much work has been dedicated to the investigation of distortion
elements of groups of diffeomorphisms of manifolds, possibly preserving
volume or symplectic forms (see for instance \cite {CF}, \cite {F},
\cite {FH}, \cite {GG}, \cite {P} and references therein),
very little seems to be known in the case of diffeomorphisms of the circle.
In this direction, Franks and Handel \cite {FH} asked whether irrational
rotations of the circle can be distortion elements.
Progress on this question had been restricted to the case of low regularity
(Calegari and Freedman \cite {CF}
give an affirmative answer in $C^1$ and remark that
their construction can not be made $C^2$).
In this note we give an affirmative answer to Franks and Handel
question in high regularity.

Let $\Diff$ be the group of orientation preserving $C^\infty$
diffeomorphisms of the circle endowed with the $C^\infty$ topology, and let
$d$ be a metric on $\Diff$ compatible with the topology.  Let us say that $f
\in \Diff$ is {\it recurrent} if $\liminf_{n \to \infty} d(f^n,\id)=0$. 
An element $f \in \Diff$ which is $C^\infty$ conjugate to a rotation is
always recurrent.  If $f$ has rational rotation number, this is the only way
to be recurrent, but for $f$ with irrational rotation number, $f$ can be
recurrent without being $C^\infty$ conjugate to rotations: indeed for every
Liouville number $\alpha$, the generic $f \in \Diff$ with rotation number
$\alpha$ is recurrent, but not $C^\infty$ conjugate
to a rotation, see \cite {Y2}.

\begin{thm}

If $f \in \Diff$ is recurrent then it is a distortion element.

\end{thm}

This is a straightforward consequence of the following estimate.

\begin{thm} \label {pr}

There exist sequences $k_n \to \infty$ and $\epsilon_n \to 0$
such that if $h_n \in \Diff$ is a sequence
such that $d(h_n,\id)<\epsilon_n$ then there exists a finitely
generated group $G$ and some $S$ such that $h_n \in G$ and
$d_S(h_n,\id) \leq k_n$ for every $n$.

\end{thm}

This implies more generally that
any countable set of recurrent elements $\{f_i\}_{i=1}^\infty \subset
\Diff$ can be made
simultaneously distortion elements of a finitely generated subgroup of
$\Diff$.\footnote {In fact, they can be made simultaneously
arbitrarily distorted, in the
sense that for every function $r:\N \to \N$, for every $i$, there exist
infinitely many $n$ such that $d_S(f_i^m,\id) \leq n$ for some
$m \geq r(n)$.  As pointed out to us by Andres Navas, in some situations
where distortion arises, it is still possible to obtain a bound on the
possible {\it distortion functions} $D_g:\N \to \N$, defined so that
$D_g(n)$ is minimal such that $\dist_S(g^m,\id)>n$ for every $m>D_g(n)$.
In particular, in the closely related setting of
piecewise-linear circle diffeomorphisms, work of Liousse-Navas shows that
distortion functions of irrational rotations grow at most exponentially.}

This kind of estimate is reminiscent of Theorem C of \cite {CF},
which is however restricted to very low regularity (homeomorphisms).
(We became familiar with \cite {CF} after the results of this paper were
obtained.  Though we have not studied their paper in detail, our approaches
seem to have considerable overlap.)

The basic strategy of this work should apply also to higher dimensions.
The only meaningful step that needs attention is the realization
of elements close to the identity as a short product of commutators of
elements close to the identity and subordinate to some fixed open covering
of the manifold (his is possibly essentially contained in the proofs of
perfectness of groups of diffeomorphisms by
Herman \cite {H}, Thurston \cite {T}, Mather \cite {M} and Epstein \cite
{E}, but we have not investigated this issue in too much detail).  On the
other hand, our approach is not appropriate to deal with volume preserving
or symplectic cases, since we use Morse-Smale dynamics in the construction.

We would like to emphasize that our work leaves open the problem of
whether irrational rotations can be distortion elements in the group of
orientation preserving
real analytic diffeomorphisms of the circle.  This question is quite
significant, and we are unwilling even to make a guess on the answer:
as pointed out to us by Andres
Navas, algebraic properties of groups of diffeomorphisms are
much more rigid (and less well understood) in the real analytic category
than in the $C^\infty$ category.

{\bf Aknowledgements:}
This work was originally motivated by a question by Rapha\"el Krikorian, who
asked the author
whether a finitely generated subgroup of $\Diff$ could be dense. 
Theorem \ref {pr} obviously addresses this question as well.
After devising the basic strategy discussed in this paper,
Etienne Ghys told me that the solution to Krikorian's question
was known (though perhaps unpublished).  Franks and Handel
question was mentioned by
Andres Navas during his course in the ICTP in
July 2008, and came to my attention through Yoccoz.  I would like to thank
Andres Navas for incentivating me to work on this problem,
for explaining some of the known theory, especially
regarding commutators, and for providing several suggestions on the writing.
This research was partially conducted
during the period the author served as Clay Research Fellow.

\section{Proof of Theorem \ref {pr}}

Let $F_1 \in \Diff$ be a map with only two
fixed points, a repelling one at $1/2$ and an attracting one at $1$, such
that $F_1$ is linear in $(1/10,9/10)$ and
outside $(1/100,99/100)$.  Let $F_2$ be a map which is the identity outside
$(1/4,3/4)$ and such that $F_2(x)=x+1/100$ for $x \in (2/5,3/5)$.  Let
$F_3$ be an irrational rotation.  We will prove the following more precise
version of Theorem \ref {pr}.

\begin{lemma} \label {better}

There exist $\epsilon_n>0$ and $k_n>0$ such that if
$h_n \in \Diff$ is such that $d(h_n,\id)<\epsilon_n$ for every
$n$ then there exist $F_4,F_5 \in \Diff$ such that
$S=\{F_1,F_2,F_3,F_4,F_5\}$
generates a group $G$ such that
$h_n \in G$ and $d_S(h_n,\id) \leq k_n$ for every $n$.

\end{lemma}

The main step in the proof of Lemma \ref {better}
is the construction of commutators
of compactly supported diffeomorphisms of an interval.

\begin{lemma} \label {commutator}

There exist $\epsilon'_n>0$ and $k'_n>0$ such that if $f_n,g_n \in
\Diff$ are
supported on $(1/7,6/7)$ and satisfy
$d(f_n,\id),d(g_n,\id)<\epsilon'_n$ for every $n$
then there exist $F_4$ and $F_5$ such that $S=\{F_1,F_2,F_3,F_4,F_5\}$
generates a group $G$
such that $[f_n,g_n] \in G$ and $d_S([f_n,g_n],\id)
\leq k'_n$.

\end{lemma}

\begin{pf}

Let $I_n$ be the interval $[1/(2n+1),1/(2n)]$, $n \geq 1$.  For each $n$,
choose $k,m \in \Z$ such that
$\tilde F_n=F_3^k F_1^m F_2 F_1^{-m} F_3^{-k}$ is
supported in $\inter I_n$ and there exists an interval
$T_n \subset I_n$ such that $\tilde F_n(x)=x+|T_n|/20$ for $x \in T_n$.

Choose $k',m' \in \Z$ such that $\hat F_n=F_3^{k'} F_1^{m'}$ takes the
interval $(1/8,7/8)$ onto an interval $T'_n \subset T_n$ of length at most
$|T_n|/40$.  Let $\Diff_a,\Diff_b \subset \Diff$ be subgroups consisting
of diffeomorphisms supported in $(1/8,7/8)$ and $T'_n$ respectively.  Then
$\Phi:\Diff_a \to \Diff_b$, $\Phi(x)=\hat F_n x \hat F^{-1}_n$ is a
continuous isomorphism.

We define $F_4$ and $F_5$ to be the identity outside the union of
the $T'_n$.  We let
$F_4=\Phi(f_n)$ and $F_5=\Phi(g_n)$.  Clearly, if $f_n$ and $g_n$ are
sufficiently close to the identity, $F_4$ and $F_5$ will indeed belong to
$\Diff$.

Consider now $A=F_4 \tilde F_n F_4^{-1} \tilde F_n^{-1}$,
$B=F_5 \tilde F_n F_5^{-1} \tilde F_n^{-1}$,
$C=F_4^{-1} F_5^{-1} \tilde F_n F_5 F_4 \tilde F_n^{-1}$.
Each of $A$, $B$ and $C$ is
supported in two intervals, $T'_n$ and $\tilde F_n(T'_n)$.  In
$T'_n$ we have $A=F_4$, $B=F_5$ and $C=F_4^{-1} F_5^{-1}$,
while in $\tilde F_n(T'_n)$ we have
$A=\tilde F_n F_4^{-1} \tilde F_n^{-1}$, $B=\tilde F_n F_5^{-1}
\tilde F_n^{-1}$, $C=\tilde F_n F_5 F_4 \tilde F_n^{-1}$.

Now $A B C$ is supported in $T'_n$,
where it is given by $F_4 F_5 F_4^{-1} F_5^{-1}$.
The map $H=\hat F^{-1}_n A B C \hat F_n$ thus satisfies $H=f_n g_n f_n^{-1}
g_n^{-1}$.  The result follows.
\end{pf}

To go from commutators to arbitrary maps, we use some explicit version
of the fact that $\Diff$ is perfect (that is, any element is a product of
commutators).


\begin{lemma} \label {commutator to map}

For every $x \in \Diff$ close to $\id$, there exists $x_1,...,x_8 \in \Diff$
close to $\id$
such that $[x_1,x_2][x_3,x_4][x_5,x_6][x_7,x_8]=x$, and
$x_1,x_2,x_7,x_8$ are supported in $(1/5,4/5)$ while
$x_3,x_4,x_5,x_6$ are supported in the
complement of $[3/10,7/10]$.

\end{lemma}

The previous lemma follows from the known theory about commutators
of diffeomorphisms (\cite {H}, \cite {T}, \cite {M}, \cite {E}),
but since the proof
is not long, we give a proof in the Appendix (in the line of Herman).

The conclusion of Lemma \ref {better} (and hence Theorem \ref {pr})
from Lemmas \ref {commutator} and \ref {commutator to map} is
straightforward.

\subsection{Proof of Lemma \ref {better}}

Let $\epsilon'_n$ and $k'_n$, $n \geq 1$ be given by Lemma \ref
{commutator}.  Choose $m>0$
such that $F_3^m(1/7,6/7)$ contains the complement of $(3/10,7/10)$.
Let $k_n=4k'_n+2m$, and let $\epsilon''_n$ be so small that if
$d(x,\id)<\epsilon''_n$ then $d(F_3^{-m} x F_3^m,\id)<\epsilon'_n$.

Using Lemma \ref {commutator to map}, choose $\epsilon_n>0$ so that
\be
\textrm {if } d(x,\id)<\epsilon_n \textrm{ then }
x=[x_1,x_2][x_3,x_4][x_5,x_6][x_7,x_8] \textrm { with }
d(x_i,\id)<\epsilon''_n
\ee
and such that $x_1,x_2,x_7,x_8$ are supported on
$(1/5,4/5)$ and $x_3,x_4,x_5,x_6$ are supported on the complement of
$[3/10,7/10]$.

Let us show that Lemma \ref {better} holds with the sequences $\epsilon_n$
and $k_n$ we have defined.

Given a sequence $h_n$ with $d(h_n,\id)<\epsilon_n$, we can write
\be \label {h_n}
h_{n+1}=[f_{4n+1},g_{4n+1}]
F_3^m [f_{4n+2},g_{4n+2}] [f_{4n+3},g_{4n+3}] F_3^{-m} [f_{4n+4},g_{4n+4}],
\ee
where $d(f_n,\id),d(g_n,\id)<\epsilon'_n$ and $f_n$
and $g_n$ are supported on $(1/7,6/7)$.  Indeed, for each
$n \geq 0$, we can just take $x=h_{n+1}$ and
$f_{4n+1}=x_1$, $g_{4n+1}=x_2$,
$f_{4n+2}=F_3^{-m} x_3 F_3^m$, $g_{4n+2}=F_3^{-m} x_4 F_3^m$,
$f_{4n+3}=F_3^{-m} x_5 F_3^m$, $g_{4n+3}=F_3^{-m} x_6 F_3^m$,
$f_{4n+4}=x_7$, $g_{4n+4}=x_8$ in (\ref {h_n}).

By Lemma \ref {commutator}, there exist $F_4$ and $F_5$ such that
$S=\{F_1,F_2,F_3,F_4,F_5\}$ generates a group $G$ with
$d_S([f_n,g_n],\id) \leq k'_n$.  It follows that $d_S(h_n,\id) \leq k_n$.

\appendix

\section{Proof of Lemma \ref {commutator to map}}

Let $\rho:\Diff \to \R/\Z$ be the rotation number.

We say that $t \in \R$ is Diophantine if $|qt-p|>C q^{-n}$ for some
$C>0$, $n>0$ and every $p \in \Z$, $q \in \Z \setminus \{0\}$.  We say that
$t \in \R/\Z$ is Diophantine if some (and thus every) lift of $t$ is
Diophantine.  We will need the following result from the Herman-Yoccoz
theory \cite {Y1}.

\begin{thm} \label {yoccoz}

If $t$ is Diophantine then for every $a \in \Diff$ such that $\rho(a)=t$,
there exists a unique $c \in \Diff$ such that $c(0)=0$ and
$c a c^{-1}$ is a rotation.  Moreover, $c$ depends continuously on $a \in
\rho^{-1}(t)$.

\end{thm}

\comm{
\begin{thm}[KAM Theorem]

Let $\gamma_t \in \Diff$, $t \in \R$ small,
be a $C^\infty$ family of diffeomorphisms such
that $\gamma_0=\id$ and $\frac {d} {dt} \gamma_t(x)<0$ for all $t$ and $x$.
Let $s \in \R$ be Diophantine.  Then
every $a \in \Diff$ near $R_s$ can be written as $\gamma_t
b R_s b^{-1}$ for some $t$ close to $0$ and some $b \in \Diff$ close
to $\id$.

\end{thm}

(The KAM Theorem is usually stated with $\gamma_t=R_{-t}$, but this
reformulation is equally easy to prove.)

\begin{thm}[Herman-Yoccoz Theorem]

If $a \in \Diff$ and $\rho(h)$ is Diophantine then there exists $c \in
\Diff$ such that $c a c^{-1}$ is a rotation.

\end{thm}

Putting together both results, we easily get:
}

It is well known that Theorem \ref {yoccoz} implies:

\begin{thm} \label {yoc}

Let $\gamma_t \in \Diff$, $t \in (-\epsilon,\epsilon)$,
be a $C^\infty$ family of diffeomorphisms such
that $\gamma_0=\id$ and $\frac {d} {dt} \gamma_t(x)<0$ for all $t$ and $x$.
If $a_0 \in \Diff$ and $\rho(h)$ is Diophantine then for every $a \in \Diff$
close to $a_0$, there exists $t \in \R$ small and $b \in \Diff$ close to
$\id$ such that $a=\gamma_t b a_0 b^{-1}$.

\end{thm}

\begin{pf}

Notice that $\rho(\gamma_t^{-1} a_0)$ and $\rho(\gamma_{-t}^{-1} a_0)$ are
small and in opposite sides of $\rho(a_0)$ for $t$ small.
Thus for every $a$ close
to $a_0$ there exists $t$ close to $0$ such that $\rho(\gamma_t^{-1}
a)=\rho(a_0)$.  The result now follows from Theorem \ref {yoccoz}.
\end{pf}

Let us consider compactly supported non-negative functions $v,w:\R/\Z \to \R$,
supported on $(1/5,4/5)$
such that $v'<0$ and $w'>0$ in $[11/50,39/50]$, $v' \leq 0$ and $w' \geq 0$
in $[21/100,79/100]$, $v=0$ on $[157/200,4/5)$ and $w=0$ on $(1/5,43/200]$. 
Letting $f_t,g_t \in \Diff$ be the flows corresponding to $-w$ and $-v$, we
see that $\frac {d} {dt} [f_t,g_t](x)>0$ for every $x \in [11/50,39/50]$ and
$\frac {d} {dt} [f_t,g_t] \geq 0$ for every $x \in \R/\Z$ and $t>0$
sufficiently small.  Let $R_{1/2}(x)=x+1/2$ be the half turn.

Now consider $h_t=[R_{1/2} f_t R_{1/2},R_{1/2} g_t R_{1/2}][f_t,g_t]$.
Then $\frac {d} {dt} h_t$ is
positive and small for $t>0$ small.  In particular, $\rho(h_t)
\neq \rho(h_0)$ for $t>0$ small (recall that $h_0=\id$).
Thus we can choose $t_0>0$
arbitrarily small such that $h_{t_0}$ has Diophantine rotation number
$\rho_0$.  Notice also that $\gamma_t=h_{t_0} h_{t_0+t}^{-1}$, $t \in \R$
small, satisfies the hypothesis of Theorem \ref {yoc}.  Applying Theorem
\ref {yoc}, we conclude
that for every $x \in \Diff$ close to $\id$, there exists $t \in \R$ small
and $b \in \Diff$ close to $\id$ such that $x=h_{t_0+t}^{-1} b h_{t_0}
b^{-1}$.  It follows that $x=[x_1,x_2][x_3,x_4][x_5,x_6][x_7,x_8]$ where
$x_1=g_{t_0+t}$, $x_2=f_{t_0+t}$, $x_3=R_{1/2} g_{t_0+t} R_{1/2}$,
$x_4=R_{1/2} f_{t_0+t} R_{1/2}$, $x_5=b R_{1/2} f_{t_0} R_{1/2} b^{-1}$,
$x_6=b R_{1/2} g_{t_0} R_{1/2} b^{-1}$, $x_7=b f_{t_0} b^{-1}$ and $x_8=b
g_{t_0} b^{-1}$, which are readily seen to be all close to the identity and
have the desired supports.

\end{document}